\pdfoutput=1
\documentclass[english]{jnsao}

\usepackage[ruled,vlined,linesnumbered,algosection,algo2e]{algorithm2e} 

\theoremstyle{definition}

\usepackage[capitalize,noabbrev,nameinlink]{cleveref}
\crefname{assumption}{Assumption}{Assumptions}
\crefname{line}{Step}{Steps}

\usepackage{enumitem}
\newlist{assumptionitemize}{enumerate}{1}
\setlist[assumptionitemize]{label = (\textsc{a}\arabic*)}
\crefname{assumptionitemizei}{Assumption}{Assumptions}
\newlist{proofitemize}{enumerate}{1}
\setlist[proofitemize]{label = (\roman*)}

\newcommand{\Rinf}{\overline{\R}}
\newcommand{\R}{\mathbb{R}}
\newcommand{\N}{\mathbb{N}}
\DeclareMathOperator*{\argmin}{\arg\min}
\newcommand{\innprod}[2]{\langle #1, #2 \rangle}
\DeclareMathOperator*{\minimize}{minimize}

\DeclareMathOperator{\dist}{dist}
\DeclareMathOperator{\stt}{subject~to}
\newcommand{\func}[3]{#1\colon#2\to#3}
\newcommand{\ffunc}[3]{#1\colon#2\rightrightarrows#3}
\DeclareMathOperator{\dom}{dom}
\DeclareMathOperator{\prox}{prox}
\DeclareMathOperator{\indicator}{\iota}
\newcommand{\XX}{\mathbb{X}}
\newcommand{\Tprox}{\operatorname{T}}
\newcommand{\Rprox}{\operatorname{R}}
\newcommand{\bigO}{\mathcal{O}}
\newcommand{\seq}[2]{\{#1\}_{#2}}
\newcommand{\cardinality}[1]{\#{#1}}

\numberwithin{equation}{section}

\newcommand{\AND}{ $\cdot$ }
\newcommand{\amsmsc}[1]{\href{http://www.ams.org/mathscinet/msc/msc2020.html?t=#1}{#1}}

\usepackage{tcolorbox}
\newtcolorbox{mybox}[1][]{%
	left=0pt,
	right=0pt,
	top=0pt,
	bottom=0pt,
	colback=gray!10,
	colframe=gray!10,
	width=\dimexpr\textwidth\relax,
	enlarge left by=0mm,
	boxsep=5pt,
	arc=5pt,outer arc=5pt,
	#1
}

\usepackage{pgfplots}
\pgfplotsset{compat=1.13}
\usepackage{tikz}
\usetikzlibrary{shapes,external}
\newcommand{\TheTitle}{Proximal gradient methods beyond monotony}
\newcommand{\TheAffiliation}{%
	University of the Bundeswehr Munich,
	Department of Aerospace Engineering,
	Institute of Applied Mathematics and Scientific Computing,
	85577 Neubiberg, Germany%
}
\newcommand{\TheKeywords}{%
	Nonsmooth nonconvex optimization\AND%
	Proximal algorithms\AND%
	Gradient methods\AND%
	Spectral stepsize\AND%
	Nonmonotone linesearch%
}
\newcommand{\TheAMSsubj}{%
	\amsmsc{49J52}\AND
	\amsmsc{65K05}\AND
 	\amsmsc{90C06}
}
\newcommand{\TheAcknowledgments}{%
	I would like to thank Axel Dreves, for valuable comments regarding the results in \cref{sec:ConvergenceRate},
	Christian Kanzow, for soliciting a clearer proof of \cref{thm:asymp}\ref{thm:asymp:cost},
	Patrick Mehlitz and Andreas Themelis, for their detailed feedback on a preliminary version of this manuscript.
	
	Part of this work was completed at Curtin University during a research visit supported by the Centre for Optimisation and Decision Science, kindly acknowledged.
	I am grateful to Ryan Loxton and Hoa Bui for their warm hospitality during my time in Perth, WA.
}

\title{\TheTitle}
\author{Alberto De Marchi%
	\thanks{\TheAffiliation. \email{alberto.demarchi@unibw.de}, \orcid{0000-0002-3545-6898}.}}
\shortauthor{De Marchi}
\manuscriptsubmitted{2022-11-10}
\manuscriptaccepted{2023-05-25}
\manuscriptvolume{4}
\manuscriptnumber{10290}
\manuscriptyear{2023}
\manuscriptdoi{10.46298/jnsao-2023-10290}

\manuscriptcopyright{{\copyright} A. De Marchi}
\manuscriptlicense{CC-BY 4.0}

\begin{document}

\maketitle

\begin{abstract}
	We address composite optimization problems, which consist in minimizing the sum of a smooth and a merely lower semicontinuous function, without any convexity assumptions.
Numerical solutions of these problems can be obtained by proximal gradient methods, which often rely on a line search procedure as globalization mechanism.
We consider an adaptive nonmonotone proximal gradient scheme based on an averaged merit function and establish asymptotic convergence guarantees under weak assumptions, delivering results on par with the monotone strategy.
Global worst-case rates for the iterates and a stationarity measure are also derived.
Finally, a numerical example indicates the potential of nonmonotonicity and spectral approximations.
	
	\medskip
	\noindent\textcolor{structure}{Keywords}\quad \TheKeywords.
	
	\medskip
	\noindent\textcolor{structure}{AMS MSC}\quad \TheAMSsubj.
\end{abstract}

\section{Introduction}\label{sec:Introduction}
In this paper we consider the classical problem of minimizing the sum of a smooth function $f$ and a nonsmooth function $g$.
More precisely, we address optimization problems of the form
\begin{equation}\label{eq:P}\tag{P}
	\minimize_{x \in \XX} \; \varphi(x),
	\qquad\qquad\text{where}\quad
	\varphi \coloneqq f + g,
\end{equation}
$\func{f}{\XX}{\R}$ is differentiable, $\func{g}{\XX}{\Rinf \coloneqq \R \cup \{\infty\} }$ admits an easily computable proximal mapping, and $\XX$ denotes an Euclidean space, i.e., a real and finite-dimensional Hilbert space.
The objective $\varphi$ and its components $f$ and $g$ are allowed to be nonconvex.

The problem class \cref{eq:P} has been widely investigated, especially in the convex setting, and a number of splitting methods have been devised for its numerical solution \cite{combettes2011proximal,parikh2014proximal,beck2017first,themelis2018forward}.
By exploiting its composite structure, proximal gradient (PG) methods offer a simple and versatile iterative technique for addressing \cref{eq:P}.
Introduced in \cite{fukushima1981generalized}, PG schemes are also known as forward-backward splittings, due to the close relation with the theory of monotone operators \cite{bauschke2017convex} and their equivalence in the convex setting \cite[\S 3]{parikh2014proximal}.
Combining a step in the negative gradient direction with a proximal point update, PG methods can cope with nonconvexity but require a mechanism to guarantee global convergence, that is, convergence from arbitrary initial points.
Widely adopted since Armijo's work \cite{armijo1966minimization} in smooth unconstrained optimization, a (backtracking) line search procedure provides such globalization, if equipped with an appropriate criterion for validating (and therefore accepting) tentative iterates.
Monotone strategies ensures that the objective value decreases \emph{at each iteration}, but possibly at the cost of taking utterly small steps.
In contrast, nonmonotone strategies adopt relaxed conditions for a tentative update to be accepted; this approach often proves beneficial in practice, as it can reduce conservatism and allow larger steps and faster convergence \cite{grippo1986nonmonotone,zhang2004nonmonotone,themelis2018forward}.
Similarly to watchdog techniques, nonmonotone line search procedures monitor progress and ensure that the objective value decreases, \emph{at least eventually}, along the iterates generated by the algorithm.

We are interested in PG methods for \cref{eq:P} with provable well-definedness and convergence guarantees under minimal assumptions on the problem data.
These properties have been recently established by Kanzow and Mehlitz \cite{kanzow2022convergence} considering an adaptive \emph{monotone} PG method.
Therein, however, the same guarantees for a \emph{nonmonotone} variant required stronger conditions, and more involved proofs.
In particular, accumulation points are shown to be stationary for \cref{eq:P} under some continuity assumptions on $\varphi$ and $g$ \cite[Ass. 4.1]{kanzow2022convergence}.

\subsection{Related Work}

In their seminal work \cite{grippo1986nonmonotone} in the context of unconstrained smooth optimization, Grippo, Lampariello and Lucidi (GLL) suggested a nonmonotone line search technique based on monitoring the \emph{maximum} objective value attained by the latest iterates.
Later Zhang and Hager \cite{zhang2004nonmonotone} proposed a different flavor of nonmonotonicity, which takes a weighted \emph{average} of the objective value at all iterates.
Anticipating the notation of \cref{alg:NMPG} below, a certain iterate $x^{k+1} \in \dom g$ is deemed acceptable upon comparison of the objective $\varphi(x^{k+1})$ with a threshold based on the merit $\Phi_k$, whose value is given by
\begin{equation}\label{eq:nm_flavours}
	\max_{j=0,\ldots,\min\{M,k\}} \varphi(x^{k-j})
	\qquad\qquad\text{and}\qquad\qquad
	(1-p) \Phi_{k-1} + p \varphi(x^k)
\end{equation}
for the max and average flavors, respectively.
Here, parameters $M \in \N$ and $p \in (0,1]$ control the level of nonmonotonicity allowed: for larger $M$, or smaller $p$,
the merit $\Phi_k$ can attain larger values, hence imposing a weaker condition for the acceptance of $x^{k+1}$.
Conversely, the values $M = 0$ and $p = 1$ result in monotone behavior, as the merit $\Phi_k$ coincides with the objective $\varphi(x^k)$.
We notice that, in general, neither of these approaches is more conservative than the other.
The GLL merit function has been adopted, e.g., in \cite{raydan1997barzilai,birgin2000nonmonotone} and \cite{wright2009sparse} for smooth and composite optimization, respectively, under some convexity assumptions, before being considered in \cite{kanzow2022convergence}.
In the context of (fully nonconvex) composite optimization, a merit function based on averaging has been adopted, e.g., in \cite{li2015accelerated} for safeguarding Nesterov-type accelerated PG methods as well as in \cite{themelis2018forward} and \cite[\S 4.3]{demarchi2022proximal} for globalizing quasi-Newton-type schemes.

Nonmonotone techniques had a significant impact also on the development and practice of spectral gradient methods. 
Long after the classical steepest descent method was proposed by Cauchy \cite{cauchy1847methode} in 1847, Barzilai and Borwein \cite{barzilai1988two} presented a simple and efficient tool for solving large-scale problems: a gradient method with step sizes motivated by Newton's method but not involving any Hessian matrix.
The Barzilai-Borwein (BB) method uses the same search direction as for Cauchy's method, but its stepsize rule is different and, at nearly no extra cost, it often outperforms standard gradient methods \cite{birgin2000nonmonotone,birgin2014spectral}.
However, even when the objective function is strongly convex, it may not yield a monotone decrease in the function value and, in fact, it may not converge at all.
Therefore, Raydan \cite{raydan1997barzilai} suggested to control this behavior by combining the spectral stepsize with a nonmonotone line search as a safeguard for globalization.
This, in particular, leads to global convergence for smooth unconstrained optimization problems and constrained on sets \cite{birgin2000nonmonotone,jia2022augmented}.

The broad class \cref{eq:P} includes constrained optimization instances, namely problems involving the minimization of a smooth cost $f(x)$ subject to the constraint $x \in X$, with $X \subseteq \XX$ some nonempty and closed set.
In fact, it suffices to take $g \equiv \indicator_X$ in \cref{eq:P}, with $\func{\indicator_X}{\XX}{\Rinf}$ denoting the indicator function of $X$, defined by $\indicator_X(x) = 0$ if $x\in X$ and $\indicator_X(x) = \infty$ otherwise.
Therefore, if the projection operator associated to $X$ (coinciding with the proximal mapping of $\indicator_X$) can be efficiently evaluated, PG methods are a viable option for dealing with such (potentially nonconvex) constrained problems.
From the literature, one can observe the (temporal) progression from BB methods \cite{barzilai1988two,raydan1997barzilai} for smooth unconstrained problems, via spectral projected gradient (SPG) methods \cite{birgin2000nonmonotone,birgin2014spectral} for smooth objectives over convex constraint sets and \cite[\S 3]{jia2022augmented} for nonconvex constraint sets, to \cite{kanzow2022convergence} for nonconvex composite objectives, all using the nonmonotone GLL line search for globalization.
Preferring a different kind of nonmonotonicity, the methods we discuss in the following allow to tackle such broad class of problems, supported by the same theoretical guarantees and under weaker assumptions.

Finally, we shall remark that these methods often play the role of inner solvers in general purpose optimization packages, where subproblems of the form \eqref{eq:P} have to be (approximately) solved.
This is the case of SPG \cite{birgin2000nonmonotone} within the augmented Lagrangian solver Algencan \cite{birgin2014practical}, \cite[Alg. 3.1]{jia2022augmented} within \cite[Alg. 4.1]{jia2022augmented}, PANOC$^{(+)}$ \cite{stella2017simple,demarchi2022proximal} within OpEn \cite{sopasakis2020open}, ALS \cite{demarchi2023constrained}, ALPS \cite{demarchi2023implicit}, and IP-FB within an interior point method \cite{demarchi2022interior}.
Thus, our findings have significant impact that possibly widens scope and aids practical performance of other methods and solvers too.

\subsection{Contribution}

Intrigued by the theoretical gap emerged in \cite{kanzow2022convergence} and motivated by the potential benefits of nonmonotone strategies in practice, we consider proximal gradient methods for addressing \cref{eq:P} in the fully nonconvex setting and focus on the role played by the line search procedure.
Why does the nonmonotone strategy lead to the difficulties highlighted in \cite[\S 4]{kanzow2022convergence}?
Is it possible to retain the simplicity and convergence guarantees of these (adaptive) methods while exploiting \emph{nonmonotone} behaviour and minimal assumptions?

As we are going to reveal, nonmonotonicity is not responsible for the theoretical issues encountered, but it is its specific flavor that may be, in some sense, too weak.
Instead of the GLL max-type condition adopted in \cite{kanzow2022convergence}, we study an adaptive PG method with a nonmonotone line search based on an \emph{averaging} merit function \cite{zhang2004nonmonotone,themelis2018forward}.
By establishing that every iteration yields some sufficient decrease for this merit function, we show that the convergence properties (and analysis) of the monotone setting naturally carry over to the nonmonotone one.

The key contributions of this paper are as follows.
\begin{itemize}
	\item We present an adaptive proximal gradient method with nonmonotone line search and establish asymptotic convergence properties and guarantees under weak assumptions on the problem data (cf. \cref{thm:NMPG:subseqStationarity}).
	Although in general not more conservative than the GLL strategy, the averaging nonmonotonicity ensures sufficient decrease at every iteration (cf. \cref{thm:NMPG:descent}), enabling a seamless extension from monotone to nonmonotone line search without additional assumptions (cf. \cref{thm:asymp:fpr}).
	\item Under local Lipschitz continuity of $\nabla f$, we provide global worst-case rate-of-convergence results.
	In particular, we show (at least) sublinear convergence of order $\bigO(1/\sqrt{k})$ for the iterates $\seq{x^k}{k\in\N}$ (cf. \cref{thm:asymp:localRateIterates}) and, along converging subsequences, for the fixed-point residual as stationarity measure (cf. \cref{thm:localLipschitzSubsequenceSublinearRate}).
	Furthermore, assuming boundedness for the set of accumulation points of $\seq{x^k}{k\in\N}$, the result is extended to the whole sequence (cf. \cref{thm:localLipschitzAccumulationCompactSublinearRate}).
\end{itemize}

\subsection*{Outline}

In the following \cref{sec:Preliminaries} we provide preliminary notions and optimality concepts.
In \cref{sec:Methods} we detail and discuss the algorithm, whose convergence properties are investigated in \cref{sec:ConvergenceAnalysis} under appropriate conditions.
We report on comparative numerical tests in \cref{sec:NumericalExample} and conclude with some remarks in \cref{sec:Conclusions}.
	
\section{Preliminaries}\label{sec:Preliminaries}
	With $\N$, $\R$ and $\Rinf \coloneqq \R \cup \{\infty\}$ we denote the natural numbers, real and extended-real line, respectively.
Throughout the paper, the Euclidean space $\XX$ will be equipped with the inner product $\func{\innprod{\cdot}{\cdot}}{\XX \times \XX}{\R}$ and the associated norm $\| \cdot \|$.
Given a nonempty set $E \subset \XX$ and a point $p\in\XX$, $\dist(p,E) \coloneqq \inf_{x\in E} \|x-p\|$ denotes the distance of $p$ from $E$.
The sum of two sets $A, B \subseteq \XX$ is meant in the sense of Minkowski, namely $A + B = B + A \coloneqq \{a+b \,\vert\, a\in A, b\in B\}$.
In case of a singleton $A = \{a\}$, we write $a + B$ as shorthand for $\{a\} + B$.

The continuous linear operator $\func{f^\prime(x)}{\XX}{\R}$ denotes the derivative of $f$ at $x \in \XX$, and we will make use of $\nabla f(x) \coloneqq f^\prime(x)^\ast 1$ where $\func{f^\prime(x)^\ast}{\R}{\XX}$ is the adjoint of $f^\prime(x)$.
This way, $\nabla f$ is a mapping from $\XX$ to $\XX$ and for any $d \in \XX$ it holds $f^\prime(x) d = \innprod{\nabla f(x)}{d}$.

The following concepts are standard in variational analysis \cite{rockafellar1998variational,mordukhovich2018variational}.
The effective \emph{domain} of an extended-real-valued function $\func{h}{\XX}{\Rinf}$ is denoted by $\dom h \coloneqq \{ x\in\XX \,\vert\, h(x) < \infty \}$.
We say that $h$ is \emph{proper} if $\dom h \neq \emptyset$
and \emph{lower semicontinuous} (lsc) if $h(\bar{x}) \leq \liminf_{x\to\bar{x}} h(x)$ for all $\bar{x} \in \XX$.
Following \cite[Def. 8.3]{rockafellar1998variational}, we denote by $\ffunc{\hat{\partial} h}{\XX}{\XX}$ the \emph{regular} (or \emph{Fr{\'e}chet}) \emph{subdifferential} of $h$, defined for $\bar{x} \in \dom h$ by
\begin{equation*}
	\hat{\partial} h(\bar{x})
	\coloneqq
	\left\{
	v \in \XX
	{}\,\middle|\,{}
	\liminf_{\substack{x\to\bar{x}\\x\neq\bar{x}}} \frac{h(x) - h(\bar{x}) - \innprod{v}{x - \bar{x}}}{\| x - \bar{x} \|} \geq 0
	\right\}
\end{equation*}
and $\hat{\partial} h(\bar{x}) \coloneqq \emptyset$ for $\bar{x} \notin \dom h$.
Then, the \emph{limiting} (or \emph{Mordukhovich}) \emph{subdifferential} of $h$ is $\ffunc{\partial h}{\XX}{\XX}$, where $v \in \partial h(\bar{x})$ if and only if $\bar{x} \in \dom h$ and there exists a sequence $\seq{x^k,v^k}{k\in\N}$ such that $x^k \in \XX$ and $v^k \in \hat{\partial} h(x^k)$ for all $k\in\N$ and $x^k \to \bar{x}$, $h(x^k) \to h(\bar{x})$, $v^k \to v$ as $k \to \infty$.
Notice that, by construction, it is always $\hat{\partial} h(x) \subset \partial h(x)$ and, whenever $h$ is convex, equality holds.
The subdifferential of $h$ at $\bar{x} \in \dom h$ satisfies $\partial(h+h_0)(\bar{x}) = \partial h(\bar{x}) + \nabla h_0(\bar{x})$ for any $\func{h_0}{\XX}{\Rinf}$ continuously differentaible around $\bar{x}$ \cite[Ex. 8.8(c)]{rockafellar1998variational}, \cite[Prop. 1.30(ii)]{mordukhovich2018variational}.

The \emph{proximal} mapping of $\func{h}{\XX}{\Rinf}$ with stepsize $\gamma>0$ is defined by
\begin{equation}
	\label{eq:proxproblem}
	\prox_{\gamma h}(x)
	{}\coloneqq{} 
	\argmin_{z\in\XX} \left\{ h(z) + \tfrac{1}{2\gamma} \| z - x \|^2 \right\} .
\end{equation}
We say that $h$ is \emph{prox-bounded} if it is proper and $h + \frac{1}{2\gamma} \| \cdot \|^2$ is bounded below on $\XX$ for some $\gamma > 0$.
The supremum of all such $\gamma$ is the threshold $\gamma_h$ of prox-boundedness for $h$.
In particular, if $h$ is bounded below by an affine function, then $\gamma_h = \infty$.
When $h$ is lsc, for any $\gamma \in (0,\gamma_h)$ the proximal mapping $\prox_{\gamma h}$ is locally bounded, nonempty- and compact-valued \cite[Thm 1.25]{rockafellar1998variational}.

	\subsection{Stationarity Concepts}\label{sec:Stationarity}
	Considering the minimization of a proper function $\func{h}{\XX}{\Rinf}$ and the definition of regular subdifferential, we have that the inclusion $0 \in \hat{\partial} h(x^\ast)$ is a necessary condition for $x^\ast \in \dom h$ to be a (local) minimizer of $h$ \cite[Thm. 10.1]{rockafellar1998variational}, \cite[Prop. 1.30(i)]{mordukhovich2018variational}.
We say a point $x^\ast \in \dom h$ is \emph{M-stationary} (for $h$) if it satisfies (the potentially weaker condition) $0 \in \partial h(x^\ast)$.
This denomination is due to the appearance of the limiting (or \emph{M}ordukhovich) subdifferential.
Moreover, a point $x^\ast\in\dom h$ is said to be \emph{$\varepsilon$-M-stationary} for some tolerance $\varepsilon \geq 0$ (or simply \emph{approximate M-stationary}) if $\dist(0, \partial h(x^\ast)) \leq \varepsilon$.
When $\varepsilon = 0$, this notion reduces to M-stationarity, as it recovers the condition $0 \in \partial h(x^\ast)$ by closedness of $\partial h(x^\ast)$ \cite[Thm. 8.6]{rockafellar1998variational}.
Notice that whenever $x^\ast \in \XX$ is an (approximate) M-stationary point relative to $h$, it necessarily belongs to its domain $\dom h$, for otherwise the subdifferential $\partial h(x^\ast)$ would be empty.

Relative to \cref{eq:P}, $\varepsilon$-M-stationarity of a point $x^\ast \in \dom \varphi$ means that it satisfies $\dist(0, \partial \varphi(x^\ast)) \leq \varepsilon$.
Under mild assumptions (see our blanket \cref{ass:f,ass:g,ass:phi} given in \cref{sec:ConvergenceAnalysis}), this condition turns into
\begin{equation}
	\dist(-\nabla f(x^\ast), \partial g(x^\ast))
	\leq
	\varepsilon
\end{equation}
because of the identity $\partial \varphi = \nabla f + \partial g$.

Finally, let us consider the proximal mapping $\prox_{\gamma h}$ defined in \cref{eq:proxproblem} and, for any given $\gamma > 0$, the proximal-gradient mapping $x \mapsto \prox_{\gamma g}(x - \gamma \nabla f(x))$ associated to \cref{eq:P}.
Then, for any $x\in \XX$, the minimizing property of $\bar{x} \in \prox_{\gamma g}(x - \gamma \nabla f(x))$ implies that $\bar{x} \in \dom g$ and
\begin{equation}\label{eq:stationarity_prox_problem}
	0
	\in
	\nabla f(x) + \frac{1}{\gamma} (\bar{x} - x) + \hat{\partial} g(\bar{x}) ,
\end{equation}
owing to the relevant necessary conditions.

\section{Methods}\label{sec:Methods}
In this section we detail and elaborate upon an adaptive nonmonotone proximal gradient scheme for solving \cref{eq:P}, whose precise statement is given in \cref{alg:NMPG}.
An investigation of its well-definedness and convergence properties under some basic assumptions is postponed to the next \cref{sec:ConvergenceAnalysis}.

Given a starting point $x^0$, the recurrence of proximal gradient updates executed at \cref{state:NMPG:barx} form the core of \cref{alg:NMPG}, whose \emph{adaptivity} stems from the backtracking line search on the ``proximal'' stepsize $\gamma_k$ at \cref{state:NMPG:ls}.
Analogously to \cite[Alg. 4.1]{kanzow2022convergence}, \cite[Alg. 2]{demarchi2022proximal}, this line search on the stepsize $\gamma_k$ obviates the need for a priori knowledge (or existence, at all) of a Lipschitz constant for $\nabla f$ and provides a globalization mechanism.
Whenever the respective needed conditions at \cref{state:NMPG:ls} are violated, the stepsize is reduced and the iteration restarted.
With the particular choice $p_k \equiv 1$ for all $k\in\N$, the merit $\Phi_k$ monitoring progress coincides with the objective $\varphi(x^k)$ for all $k\in\N$, thus requiring the objective value to monotonically decrease along the iterates due to \cref{state:NMPG:ls}.
Instead, the sequence $\seq{ \varphi(x^k) }{k\in\N}$ can exhibit some degree of \emph{nonmonotonicity} whenever $p_k < 1$, despite the monotone decrease of the monitoring sequence $\seq{ \Phi_k }{k\in\N}$ generated by the \emph{averaging} at \cref{state:NMPG:Phi}.

\begin{algorithm2e}[tbh]
	\DontPrintSemicolon
	\caption{Proximal Gradient Method for \cref{eq:P}}%
	\label{alg:NMPG}%
	\KwData{starting point $x^0 \in \dom g$, termination tolerance $\varepsilon > 0$}
	\KwResult{$\varepsilon$-stationary point $x^\ast$ for \cref{eq:P}}
	select $0 < \gamma_{\min} \leq \gamma_{\max} < \infty$,
		$0 < \alpha_{\min} \leq \alpha_{\max} < 1$,
		$0 < \beta_{\min} \leq \beta_{\max} < 1$,
		$p_{\min} \in (0,1]$\;
	set $\Phi_0 \gets \varphi(x^0)$\;
	\For{$k = 1,2,\ldots$}{
		choose $\gamma_k \in [\gamma_{\min}, \gamma_{\max}]$\label{state:NMPG:init}\;
		compute $x^k \in \prox_{\gamma_k g}( x^{k-1} - \gamma_k \nabla f(x^{k-1}) )$\label{state:NMPG:barx}\;
		\If{$\| \frac{1}{\gamma_k} (x^k - x^{k-1}) - \nabla f(x^k) + \nabla f(x^{k-1}) \| \leq \varepsilon$\label{state:NMPG:termination}}{
			\Return $x^\ast \gets x^k$
		}%
		choose $\alpha_k \in [\alpha_{\min}, \alpha_{\max}]$ and $\beta_k \in [\beta_{\min}, \beta_{\max}]$\label{state:NMPG:alpha}\label{state:NMPG:beta}\;
		\If{$\varphi(x^k) > \Phi_{k-1} - \frac{1 - \alpha_k}{2 \gamma_k} \| x^k - x^{k-1} \|^2$\label{state:NMPG:ls}}{
			set $\gamma_k \gets \beta_k \gamma_k$ and go back to \cref{state:NMPG:barx}\label{state:NMPG:gamma}
		}%
		choose $p_k \in [p_{\min}, 1]$ and set $\Phi_k \gets (1 - p_k) \Phi_{k-1} + p_k \varphi(x^k)$\label{state:NMPG:Phi}\;
	}
\end{algorithm2e}

The procedure relies on many algorithmic parameters $\gamma_k$, $\alpha_k$, $\beta_k$ and $p_k$, which can be selected anew (possibly several times) at every iteration $k$, at \cref{state:NMPG:init,state:NMPG:alpha,state:NMPG:Phi}.
These degrees of freedom provide significant flexibility in controlling the behavior of \cref{alg:NMPG} and the iterates it generates, and possibly in steering and speeding up the convergence.

It should be highlighted that it is not restrictive to require a feasible starting point for \cref{alg:NMPG}, namely a point $x^0 \in \dom g$.
Conversely, it simplifies exposition and discussion of the convergence analysis in the following \cref{sec:ConvergenceAnalysis}.
That $x^0 \in \dom g$ comes without loss of generality follows from the observation that starting with some $x^0 \notin \dom g$ would result in a trivial iteration.
More precisely, \cref{state:NMPG:barx} would generate some $x^1 \in \dom g$ that is accepted by the line search at \cref{state:NMPG:ls} without need for backtracking, since $\varphi(x^1) < \Phi_0 = \varphi(x^0) = \infty$.
However, in order to avoid the infinite-valued merit $\Phi_0$ spoiling the line search condition (due to \cref{state:NMPG:Phi}, it would be $\Phi_k = \infty$ for all $k\in\N$), the procedure should be restarted with the newly set $x^0 \gets x^1$.

\subsection{Termination Criteria}

\cref{state:NMPG:termination} of \cref{alg:NMPG} entails the termination criterion, which is designed to detect approximate stationary points for \cref{eq:P}.
Rearranging the necessary condition \cref{eq:stationarity_prox_problem} associated to the proximal mapping evaluated at \cref{state:NMPG:barx}, we observe that the iterates satisfy
\begin{equation*}
	\frac{1}{\gamma_k} (x^{k-1} - x^k) - \nabla f(x^{k-1}) + \nabla f(x^k)
	\in 
	\nabla f(x^k) + \hat{\partial} g(x^k)
	=
	\hat{\partial} \varphi(x^k)
\end{equation*}
for all $k\in\N$.
This inclusion justifies the termination condition at \cref{state:NMPG:termination} since
\begin{equation}\label{eq:boundDistanceSubdifferential}
	\dist(0,\partial \varphi(x^k))
	\leq
	\dist(0,\hat{\partial}\varphi(x^k))
	\leq
	\left\| \frac{1}{\gamma_k} (x^{k-1} - x^k) - \nabla f(x^{k-1}) + \nabla f(x^k) \right\|
	\leq
	\varepsilon
\end{equation}
encodes approximate stationarity of $x^k$ for \cref{eq:P}.
Notice that the first inequality in \cref{eq:boundDistanceSubdifferential} follows from the inclusion $\hat{\partial} \varphi(x^k) \subset \partial \varphi(x^k)$.
As established in \cref{thm:finiteTermination} below, for any given $\varepsilon > 0$ the condition at \cref{state:NMPG:termination} is satisfied for some $k \in \N$ large enough, thus allowing the algorithm to return with the corresponding $\varepsilon$-M-stationary point $x^k \in \dom \varphi$.

Finally, we shall point out that the termination condition is checked within the backtracking loop, as opposed to after the line search procedure; analogously to \cite[Alg. 3.1]{jia2022augmented} and \cite[Alg. 4.1]{kanzow2022convergence}, respectively.
Although potentially costly, this allows to avoid a potential infinite loop at \cref{state:NMPG:ls}, thus making \cref{alg:NMPG} well-defined.
If the problem data in \cref{eq:P} satisfy stronger regularity properties (e.g., $\nabla f$ is locally Lipschitz continuous, cf. \cref{thm:finite:LS}), it may be possible to safely move the termination checking out of the backtracking procedure, thereby requiring fewer function evaluations.
	
\section{Convergence Analysis}\label{sec:ConvergenceAnalysis}
	In this section we analyze the properties of the iterates generated by \cref{alg:NMPG}.
The following blanket assumptions on \cref{eq:P} are considered throughout.
\begin{mybox}
	\begin{assumptionitemize}
		\item\label{ass:f}%
		$\func{f}{\XX}{\R}$ is continuously differentiable;
		\item\label{ass:g}%
		$\func{g}{\XX}{\Rinf}$ is proper, lsc and prox-bounded;
		\item\label{ass:phi}%
		$\inf \varphi \in \R$.
	\end{assumptionitemize}
\end{mybox}
Notice that \cref{ass:phi} and properness and lsc of $g$ in \cref{ass:g} are relevant for the well-definedness of \cref{eq:P}, as they imply that $\varphi \coloneqq f + g$ is lsc and bounded from below on $\dom \varphi = \dom g \ne \emptyset$.
\Cref{ass:f} and prox-boundedness of $g$ in \cref{ass:g} guarantee that the proximal mapping evaluations and the algorithm are well-defined, as we are about to show.
Furthermore, since the iterates generated by \cref{alg:NMPG} satisfy $\seq{x^k}{k\in\N} \subseteq \dom g$, the smooth function $f$ could be defined just on the domain of $g$.

The first basic result establishes that the algorithm is well-defined, namely that the inner (line search) loop requires finitely many backtrackings.
In \cref{thm:finite:LS}\ref{thm:finite:LS:nonMstationary}--\ref{thm:finite:LS:locallyLipschitz} the finite termination of the line search procedure is due to the tentative update yielding sufficient decrease according to the condition at \cref{state:NMPG:ls}, whereas in \cref{thm:finite:LS}\ref{thm:finite:LS:epsnonzero} it follows from \cref{state:NMPG:termination} detecting that the current iterate is $\varepsilon$-M-stationary for \cref{eq:P}.
\begin{mybox}
	\begin{lemma}[Well-definedness]\label{thm:finite:LS}
		Suppose that \crefrange{ass:f}{ass:phi} are satisfied.		
		Consider the $k$-th iteration of \cref{alg:NMPG} and assume that at least one of the following conditions holds:
		\begin{proofitemize}
			\item\label{thm:finite:LS:nonMstationary} $x^{k-1}$ is not M-stationary for \cref{eq:P};
			\item\label{thm:finite:LS:locallyLipschitz} $\nabla f$ is locally Lipschitz continuous (at least in a neighborhood of $x^{k-1}$);
			\item\label{thm:finite:LS:epsnonzero} $\varepsilon > 0$.
		\end{proofitemize}
		Then, the iteration terminates, and in particular \cref{state:NMPG:ls} is passed in finitely many backtrackings.
	\end{lemma}
\end{mybox}
\begin{proof}
	The line of proof builds upon \cite[Lemma 3.1]{kanzow2022convergence}, relevant for the monotone counterpart of \cref{alg:NMPG}, and the properties of the monitoring sequence $\seq{\Phi_k}{k\in\N}$ with respect to $\seq{\varphi(x^k)}{k\in\N}$.
	\begin{itemize}
		\item[\ref{thm:finite:LS:nonMstationary}]
		The claim follows immediately by induction since $\Phi_0 \coloneqq \varphi(x^0)$ and, for $k \geq 1$, the conditions at \cref{state:NMPG:Phi,state:NMPG:ls} imply
		\begin{align}\label{eq:phi_leq_Phi}
			\Phi_k
			{}\coloneqq{}
			(1 - p_k) \Phi_{k-1} + p_k \varphi(x^k)
			{}\geq{}&
			(1 - p_k) \left( \varphi(x^k) + \frac{1 - \alpha_k}{2 \gamma_k} \| x^k - x^{k-1} \|^2 \right) + p_k \varphi(x^k) \\
			{}={}&
			\varphi(x^k) + (1 - p_k) \frac{1 - \alpha_k}{2 \gamma_k} \| x^k - x^{k-1} \|^2
			\geq
			\varphi(x^k) , \nonumber
		\end{align}
		proving that the nonmonotone line search is not more conservative than the monotone one (which has $p_k = 1$ and $\Phi_k = \varphi(x^k)$ for all $k\in\N$).
		Thus, relying on \cite[Lemma 3.1]{kanzow2022convergence}, the condition at \cref{state:NMPG:ls} is violated in finitely many steps.
		\item[\ref{thm:finite:LS:locallyLipschitz}]
		Since $x^{k-1}$ remains untouched during the $k$-th iteration, $\gamma_k$ is upper bounded by $\gamma_{\max}$, and $\ffunc{\prox_{\gamma_k g}}{\XX}{\XX}$ is locally bounded, we have that the (tentative) updates $x^k$ at \cref{state:NMPG:barx} remain in a bounded set \cite[Lemma 4.1]{demarchi2022proximal}; let $L_k$ be the (local) Lipschitz constant of $\nabla f$ on this set.
		Therefore, for any given $\alpha_k \in (0,1)$ the quadratic upper bound
		\begin{equation}\label{eq:finite:quadupperbound}
			f(x^k) \leq f(x^{k-1}) + \innprod{\nabla f(x^{k-1})}{x^k - x^{k-1}} + \frac{\alpha_k}{2 \gamma_k} \| x^k - x^{k-1} \|^2
		\end{equation}
		holds after finitely many backtrackings, at most when $\gamma_k \leq \alpha_k / L_k$, by \cite[Lemma 5.7]{beck2017first}.
		Noticing that the minimizing property of $x^k$ at \cref{state:NMPG:barx} gives
		\begin{equation*}
			g (x^k) + \innprod{\nabla f(x^{k-1})}{x^k - x^{k-1}} + \frac{1}{2 \gamma_k} \| x^k - x^{k-1} \|^2 \leq g (x^{k-1}) ,
		\end{equation*}
		adding terms with \cref{eq:finite:quadupperbound} leads to
		\begin{equation*}
			\varphi(x^k)
			{}\leq{}
			\varphi(x^{k-1}) - \frac{1 - \alpha_k}{2 \gamma_k} \| x^k - x^{k-1} \|^2
			{}\leq{}
			\Phi_{k-1} - \frac{1 - \alpha_k}{2 \gamma_k} \| x^k - x^{k-1} \|^2
			,
		\end{equation*}
		where the last inequality follows from the upper bound \cref{eq:phi_leq_Phi}.
		This proves that the condition at \cref{state:NMPG:ls} is violated and the inner loop terminates.
		\item[\ref{thm:finite:LS:epsnonzero}]
		In view of \cref{thm:finite:LS}\ref{thm:finite:LS:nonMstationary}, let us assume that $x^{k-1}$ is M-stationary.
		Then, by \cref{ass:f} and $\varepsilon > 0$, the termination condition at \cref{state:NMPG:termination} is satisfied after finitely many attempts owing to the arguments of \cite[Lemma 3.1]{kanzow2022convergence}.
		\qedhere
	\end{itemize}
\end{proof}
We already mentioned that $\seq{ \varphi(x^k) }{k\in\N}$ may possess a nonmonotone behavior, whereas the sequence $\seq{ \Phi_k }{k\in\N}$ monitoring (and forcing) progress along the iterates is monotonically decreasing.
As a corollary, the iterates of \cref{alg:NMPG} belong to the sublevel set of $\varphi$ induced by the starting point $x^0 \in \dom g$.
Furthermore, we establish that $\seq{ \Phi_k }{k\in\N}$ exhibits some sufficient decrease at every iteration.
This feature appears to provide the foundation for proving convergence under weak assumptions also in the nonmonotone setting, thus bridging the gap emerged in \cite{kanzow2022convergence}.
\begin{mybox}
	\begin{lemma}[Descent behavior]\label{thm:NMPG:descent}
		Suppose that \crefrange{ass:f}{ass:phi} are satisfied.
		The following hold for the sequence of iterates $\seq{ x^k }{k\in\N}$ generated by \cref{alg:NMPG}:
		\begin{proofitemize}
			\item\label{thm:descent:descent}%
			The sequence $\seq{ \Phi_k }{k\in\N}$ is monotonically decreasing, lower bounded by $\inf\varphi$ and, for every $k \geq 1$, $k \in \N$, one has
			\begin{equation}\label{eq:SD}
				\varphi(x^k) + (1 - p_k) \delta_k
				\leq
				\Phi_k
				\leq
				\Phi_{k-1} - p_k \delta_k
				\qquad\text{where}\quad
				\delta_k
				{}\coloneqq{}
				\frac{1 - \alpha_k}{2 \gamma_k}
				\| x^k - x^{k-1} \|^2 .
			\end{equation}
			\item\label{thm:descent:sublevel}%
			Every iterate $x^k$ remains in the sublevel set $\{ x\in\XX \,\vert\, \varphi(x)\leq\varphi(x^0) \} \subseteq \dom g$.
		\end{proofitemize}
	\end{lemma}
\end{mybox}
\begin{proof}
	Based on \cref{thm:finite:LS}, we may assume the iterates of \cref{alg:NMPG} are well-defined.
	\begin{enumerate}
		\item[\ref{thm:descent:descent}]
		By the update rule at \cref{state:NMPG:Phi} and the line search condition at \cref{state:NMPG:ls}, we have that
		\begin{align*}
			\Phi_k
			{}\coloneqq{}&
			(1-p_k) \Phi_{k-1} + p_k \varphi(x^k)
			\\
			{}\leq{}&
			(1-p_k) \Phi_{k-1} + p_k \left( \Phi_{k-1} - \frac{1 - \alpha_k}{2 \gamma_k} \| x^k - x^{k-1} \|^2 \right) \\
			{}={}&
			\Phi_{k-1} - p_k \frac{1 - \alpha_k}{2 \gamma_k} \| x^k - x^{k-1} \|^2 .
		\end{align*}
		The lower bound on $\Phi_k$ follows directly from \cref{ass:phi} and the inequalities in \cref{eq:phi_leq_Phi}.
		\item[\ref{thm:descent:sublevel}]
		Owing to $x^0 \in \dom g$ and continuity of $f$, the sublevel set of $\varphi$ induced by $x^0$ is a subset of $\dom g$.
		Then, since $p_k \in (0,1]$ and $\delta_k \geq 0$, it follows from $\Phi_0 \coloneqq \varphi(x^0)$ and the inequality in \cref{eq:SD} that $\varphi(x^k)
		\leq
		\Phi_k
		\leq\ldots\leq
		\Phi_0
		=
		\varphi(x^0)
		<
		\infty$,
		whence the inclusion of $x^k$ in the sublevel set.
		\qedhere
	\end{enumerate}
\end{proof}
We next consider an asymptotic analysis of the algorithm.
The following \cref{thm:asymp}\ref{thm:asymp:diff_x} is the counterpart of \cite[Prop 3.1]{kanzow2022convergence} in the nonmonotone setting.
Contrarily to \cite[Prop 4.1]{kanzow2022convergence}, it is not more difficult to prove, thanks to the sufficient decrease on the averaging merit function; cf. \cref{eq:SD}.
\begin{mybox}
	\begin{lemma}\label{thm:asymp}%
		Suppose that \crefrange{ass:f}{ass:phi} are satisfied.
		The following hold for the sequence of iterates $\seq{ x^k }{k\in\N}$ generated by \cref{alg:NMPG}:
		\begin{proofitemize}
			\item\label{thm:asymp:cost}%
			$\seq{ \varphi(x^k) }{k\in\N}$ and $\seq{ \Phi_k }{k\in\N}$ converge to a finite value $\varphi_\star \geq \inf \varphi$; the latter from above.
			\item\label{thm:asymp:summable}%
			$
			\sum_{k\in\N} \frac{1}{\gamma_k} \| x^k - x^{k-1} \|^2
			{}<{}
			\infty
			$.
			\item\label{thm:asymp:diff_x}%
			$
			\lim_{k\to\infty} \| x^k - x^{k-1} \|
			{}={}
			0
			$.
		\end{proofitemize}
	\end{lemma}
\end{mybox}
\begin{proof}
	We implicitly assume that \cref{alg:NMPG} generates an infinite sequence $\seq{ x^k }{k\in\N}$. 
	\begin{enumerate}
		\item[\ref{thm:asymp:cost}]%
		The assertion regarding $\seq{\Phi_k}{k\in\N}$ follows from \cref{eq:SD} and \cref{ass:phi}.
		Then, consider the update rule at \cref{state:NMPG:Phi}, which gives $\varphi(x^k) = \tfrac{1}{p_k} (\Phi_k - \Phi_{k-1}) + \Phi_{k-1}$.
		Taking the limit for $k\to\infty$ and noticing that $p_k\geq p_{\min}>0$ for all $k\in\N$, the convergence of $\seq{\Phi_k}{k\in\N}$ implies that of $\seq{\varphi(x^k)}{k\in\N}$ to the same finite value $\varphi_\star$.
		\item[\ref{thm:asymp:summable}]%
		A telescoping argument on \cref{eq:SD}, together with \ref{thm:asymp:cost}, yields
		\begin{equation}\label{eq:NMPG:asymp:summable}
			\Phi_0 - \varphi_\star
			\geq
			\sum_{j=1}^k \Phi_{j-1} - \Phi_j
			\geq
			\sum_{j=1}^k p_j \delta_j
			\geq
			\sum_{j=1}^k p_{\min} \frac{1 - \alpha_{\max}}{2 \gamma_j} \| x^j - x^{j-1} \|^2 ,
		\end{equation}
		since $p_j \geq p_{\min}$ and $\alpha_j \leq \alpha_{\max}$ for all $j \geq 1$.
		The claimed finite sum follows from $p_{\min} > 0$, $\alpha_{\max} < 1$, and the independence of the finite upper bound $\Phi_0 - \varphi_\star$ on $k \in \N$.
		\item[\ref{thm:asymp:diff_x}]%
		Follows from \ref{thm:asymp:summable}, since $\seq{\gamma_k}{k\in\N}$ is upper bounded by $\gamma_{\max}$.
		\qedhere
	\end{enumerate}
\end{proof}
The following result is the nonmonotone counterpart of \cite[Prop 3.2]{kanzow2022convergence} and provides a key ingredient in relation to (the detection of) M-stationarity; see the termination condition at \cref{state:NMPG:termination}.
In stark contrast with \cite[Thm 4.1, Ass. 4.1]{kanzow2022convergence}, it does not require additional continuity assumptions on $\varphi$ nor $g$, and its proof closely patterns that of \cite[Prop 3.2]{kanzow2022convergence}, with minor modifications to account for the nonmonotone line search.
\begin{mybox}
	\begin{lemma}\label{thm:asymp:fpr}
		Suppose that \crefrange{ass:f}{ass:phi} are satisfied.
		Consider a sequence $\seq{x^k}{k\in\N}$ generated by \cref{alg:NMPG}.
		Let $x^\star$ be an accumulation point of $\seq{ x^k }{k\in\N}$ and $\seq{ x^k }{k\in K}$ a subsequence such that $x^k \to_K x^\star$.
		Then, it is $\frac{1}{\gamma_k} \| x^k - x^{k-1} \| \to_K 0$.
	\end{lemma}
\end{mybox}
\begin{proof}
	If the subsequence $\seq{ \gamma_k }{k\in K}$ is bounded away from zero, the statement follows immediately from \cref{thm:asymp}\ref{thm:asymp:diff_x}.
	The remaining part of this proof therefore assumes, without loss of generality, that $\gamma_k  \to_K 0$.
	Closely following the proof of \cite[Prop. 3.2]{kanzow2022convergence}, let $\hat{\gamma}_k$ denote the attempted stepsize preceding the accepted stepsize $\gamma_k$; these values are well-defined as the inner loop terminates by \cref{thm:finite:LS}.
	Owing to \cref{state:NMPG:gamma}, it is $\gamma_k / \beta_{\max} \leq \hat{\gamma}_k \leq \gamma_k / \beta_{\min}$ for each $k\in\N$ and so, by $\gamma_k  \to_K 0$, we also have $\hat{\gamma}_k \to_K 0$.
	However, the corresponding vector $\hat{x}^k$ does not satisfy the stepsize condition from \cref{state:NMPG:ls} (for otherwise $\hat{\gamma}_k$ would have been accepted) with the associated parameter $\hat{\alpha}_k \in [\alpha_{\min}, \alpha_{\max}]$, namely it is
	\begin{equation}\label{Prop:gammaxdiff_NM-1}
		\varphi (\hat{x}^k)
		>
		\Phi_{k-1} - \frac{1 - \hat{\alpha}_k}{2 \hat{\gamma}_k} \| \hat{x}^k - x^{k-1} \|^2
		\geq
		\Phi_{k-1} - \frac{1 - \alpha_{\min}}{2 \hat{\gamma}_k} \| \hat{x}^k - x^{k-1} \|^2 .
	\end{equation}
	On the other hand, since $ \hat{x}^k $ solves the corresponding subproblem at \cref{state:NMPG:barx}, we have
	\begin{equation}\label{Prop:gammaxdiff_NM-2}
		\innprod{\nabla f(x^{k-1})}{\hat{x}^k - x^{k-1}} + \frac{1}{2 \hat{\gamma}_k} \| \hat{x}^k - x^{k-1} \|^2 + g( \hat{x}^k )
		\leq
		g( x^{k-1} ) .
	\end{equation}
	Then, using the Cauchy-Schwarz inequality and \cref{thm:NMPG:descent}\ref{thm:descent:sublevel}, we obtain
	\begin{align*}
		g( \hat{x}^k ) + \frac{1}{2 \hat{\gamma}_k} \| \hat{x}^k - x^{k-1} \|^2 
		{}\leq{}& 
		\| \nabla f(x^{k-1}) \| \| \hat{x}^k - x^{k-1} \| + g(x^{k-1}) \\
		{}={}&
		\| \nabla f(x^{k-1}) \| \| \hat{x}^k - x^{k-1} \| + \varphi(x^{k-1}) - f(x^{k-1}) \\
		{}\leq{}& 
		\| \nabla f(x^{k-1}) \| \| \hat{x}^k - x^{k-1} \| + \Phi_0 - f(x^{k-1}) .
	\end{align*}
	By continuous differentiability of $f$ and prox-boundedness of $g$, if $\seq{ \| \hat{x}^k - x^{k-1} \| }{k \in K}$ is unbounded, then the left-hand side grows more rapidly than the right-hand side.
	Analogously, if $\seq{ \| \hat{x}^k - x^{k-1} \| }{k\in K}$ remains bounded, but staying away from zero, at least on a subsequence, then the right-hand side is bounded but the left-hand side is unbounded, on the corresponding subsequence, as $\hat{\gamma}_k \to_K 0$.
	Therefore, it must be that $\hat{x}^k - x^{k-1} \to_K 0$, and thus $\hat{x}^k \to_K x^\star$ by \cref{thm:asymp}\ref{thm:asymp:diff_x}.

	Now, by the mean-value theorem, there exists $\xi^k$ on the line segment connecting $x^{k-1}$ with $\hat{x}^k$ such that
	\begin{align*}
		\varphi( \hat{x}^k ) - \varphi(x^{k-1}) 
		{}={}& 
		f( \hat{x}^k ) + g(\hat{x}^k) - f(x^{k-1}) - g(x^{k-1}) \\
		{}={}&
		\innprod{ \nabla f(\xi^k)}{\hat{x}^k - x^{k-1}} + g(\hat{x}^k) - g(x^{k-1}) .
	\end{align*}
	Substituting $g(\hat{x}^k) - g(x^{k-1})$ from this expression into \cref{Prop:gammaxdiff_NM-2} yields 
	\begin{equation*}
		\innprod{\nabla f(x^{k-1}) - \nabla f( \xi^k)}{\hat{x}^k - x^{k-1}} + \frac{1}{2 \hat{\gamma}_k} \| \hat{x}^k - x^{k-1} \|^2 + \varphi( \hat{x}^k )
		\leq
		\varphi(x^{k-1}) .
	\end{equation*}
	Therefore, exploiting \cref{Prop:gammaxdiff_NM-1} and \cref{thm:NMPG:descent}\ref{thm:descent:descent}, we obtain
	\begin{align*}
		\innprod{\nabla f(x^{k-1}) - \nabla f( \xi^k)}{\hat{x}^k - x^{k-1}}
		{}+{}
		\frac{1}{2 \hat{\gamma}_k} \| \hat{x}^k - x^{k-1} \|^2
		{}\leq{}&
		\varphi(x^{k-1}) - \varphi( \hat{x}^k ) \\
		{}\leq{}&
		\varphi(x^{k-1}) - \Phi_{k-1} + \frac{1 - \alpha_{\min}}{2 \hat{\gamma}_k} \| \hat{x}^k - x^{k-1} \|^2 \\
		{}\leq{}&
		\frac{1 - \alpha_{\min}}{2 \hat{\gamma}_k} \| \hat{x}^k - x^{k-1} \|^2 .
	\end{align*}
	Then, after rearranging, the Cauchy-Schwarz inequality yields
	\begin{equation*}
		\frac{\alpha_{\min}}{2 \hat{\gamma}_k} \| \hat{x}^k - x^{k-1} \|^2
		{}\leq{}
		\| \nabla f(x^{k-1}) - \nabla f(\xi^k) \| \| \hat{x}^k - x^{k-1} \| .
	\end{equation*}
	Since $\hat{x}^k \neq x^{k-1}$ in view of \cref{Prop:gammaxdiff_NM-1} and \cref{thm:NMPG:descent}\ref{thm:descent:descent}, the previous inequality implies that
	\begin{equation}\label{Prop:gammaxdiff_NM-4}
		\frac{\alpha_{\min}}{2 \hat{\gamma}_k} \| \hat{x}^k - x^{k-1} \|
		\leq
		\| \nabla f(x^{k-1}) - \nabla f(\xi^k) \| .
	\end{equation}
	Owing to $x^k\to_K x^\star$ and $\hat{x}^k \to_K x^\star$, it must be also $\xi^k \to_K x^\star$.
	Therefore, using the continuous differentiability of $f$ and the fact that $\alpha_{\min} > 0$, it follows from \cref{Prop:gammaxdiff_NM-4} that $\frac{1}{\hat{\gamma}_k} \| \hat{x}^k - x^{k-1} \| \to_K 0$.

	We now move toward the assertion.
	Exploiting the fact that $x^k$ and $\hat{x}^k$ are solutions of the subproblems at \cref{state:NMPG:barx} with stepsize $\gamma_k$ and $\hat{\gamma}_k$, respectively, we obtain that
	\begin{equation*}
		\innprod{\nabla f(x^{k-1})}{x^k - \hat{x}^k} + \frac{1}{2 \gamma_k} \| x^k - x^{k-1} \|^2 + g(x^k)
		{}\leq{}
		\frac{1}{2 \gamma_k} \| \hat{x}^k - x^{k-1} \|^2 + g(\hat{x}^k)
	\end{equation*}
	and
	\begin{equation*}
		\innprod{\nabla f(x^{k-1})}{\hat{x}^k - x^k} + \frac{1}{2 \hat{\gamma}_k} \| \hat{x}^k - x^{k-1} \|^2 + g(\hat{x}^k)
		{}\leq{}
		\frac{1}{2 \hat{\gamma}_k} \| x^k - x^{k-1} \|^2 + g(x^k) .
	\end{equation*}
	Owing to $\gamma_k \leq \beta_{\max} \hat{\gamma}_k < \hat{\gamma}_k$, adding these two inequalities yields the bound $\| x^k - x^{k-1} \| \leq \| \hat{x}^k - x^{k-1} \|$.
	Therefore, by $\gamma_k \geq \beta_{\min} \hat{\gamma}_k > 0$, we have that
	\begin{equation*}
		\frac{1}{\gamma_k} \| x^k - x^{k-1} \|
		\leq
		\frac{1}{\beta_{\min} \hat{\gamma}_k} \| x^k - x^{k-1} \|
		\leq
		\frac{1}{\beta_{\min} \hat{\gamma}_k} \| \hat{x}^k - x^{k-1} \|
		\to_K 0 ,
	\end{equation*}
	completing the proof.
\end{proof}
Patterning the proof above, local Lipschitz continuity of $\nabla f$ implies that the stepsize $\gamma_k$ does not vanish along a convergent subsequence.
This observation is stated explicitly in the  following result, analogous to \cite[Cor. 3.1]{kanzow2022convergence}.
A detailed proof is included in the Additional Material (p.\ \pageref{proof:cor:gammaLocalLipschitz}).
\begin{mybox}
	\begin{corollary}\label{cor:gammaLocalLipschitz}
		Suppose that \crefrange{ass:f}{ass:phi} are satisfied and consider a sequence $\seq{x^k}{k\in\N}$ generated by \cref{alg:NMPG}.
		Let $x^\star$ be an accumulation point of $\seq{x^k}{k\in\N}$ and $\seq{x^k}{k\in K}$ a subsequence such that $x^k \to_K x^\star$.
		If $\nabla f$ is locally Lipschitz continuous (at least in a neighborhood of $x^\star$), then the subsequence $\seq{\gamma_k}{k\in K}$ is bounded away from zero.
	\end{corollary}
\end{mybox}
The following result provides the main global convergence guarantees for \cref{alg:NMPG}.
Notice that additional smoothness conditions are required on either $f$ or $g$.
These match those considered in \cite[Thm 3.1]{kanzow2022convergence} for the monotone variant.
\begin{mybox}
	\begin{theorem}\label{thm:NMPG:subseqStationarity}
		Suppose that \crefrange{ass:f}{ass:phi} hold and at least one of the following conditions is satisfied for \cref{eq:P}:
		\begin{proofitemize}
			\item\label{ass:subseqStationarity:gContinuous}%
			$\func{g}{\XX}{\Rinf}$ is continuous relative to its domain;
			\item\label{ass:subseqStationarity:fLocallyLipschitz}%
			$\func{\nabla f}{\XX}{\XX}$ is locally Lipschitz continuous.
		\end{proofitemize}
		Consider a sequence $\seq{ x^k }{k\in\N}$ generated by \cref{alg:NMPG} and the finite value $\varphi_\star$ as in \cref{thm:asymp}\ref{thm:asymp:cost}.
		Then, each accumulation point $x^\star$ of $\seq{ x^k }{k\in\N}$ is M-stationary for \cref{eq:P} and it is $\varphi(x^\star) = \varphi_\star$.
	\end{theorem}
\end{mybox}
\begin{proof}
	Let $\seq{ x^k }{k\in K}$ be a subsequence converging to $x^\star$.
	In view of \cref{thm:asymp}\ref{thm:asymp:diff_x}, it follows that also the subsequence $\seq{ x^{k-1} }{k\in K}$ converges to $x^\star$.
	Furthermore, \cref{thm:asymp:fpr} yields $\frac{1}{\gamma_k} \| x^k - x^{k-1} \| \to_K 0$.
	Hence, if we can show the $g$-attentive convergence of $\seq{x^k}{k\in K}$ to $x^\star$, namely $x^k \to_K x^\star$ with $g(x^k) \to_K g(x^\star)$, the desired statement $0 \in \nabla f(x^\star) + \partial g(x^\star)$ is due to the optimality conditions \cref{eq:stationarity_prox_problem} for \cref{state:NMPG:barx} and is obtained by taking the limit $k \to_K \infty$ in
	\begin{equation}\label{eq:stationarity_prox_problem_indexed}
		0
		\in
		\nabla f(x^{k-1}) + \frac{1}{\gamma_k} (x^k - x^{k-1}) + \hat{\partial} g(x^k) .
	\end{equation}
	Moreover, by continuity of $f$ and \cref{thm:asymp}\ref{thm:asymp:cost}, it also implies that $\varphi(x^k) \to_K \varphi(x^\star) = \varphi_\star$.
	\begin{enumerate}
		\item[\ref{ass:subseqStationarity:gContinuous}]%
		If $g$ is continuous on its domain, the condition $g(x^k) \to_K g(x^\star)$ is readily obtained since all iterates $x^k$ generated by \cref{alg:NMPG} belong to $\dom g$ and, owing to \cref{thm:NMPG:descent}\ref{thm:descent:sublevel} and lsc of $g$ by \cref{ass:g}, $x^\star \in \dom \varphi = \dom g$.
		\item[\ref{ass:subseqStationarity:fLocallyLipschitz}]%
		It remains to consider the situation where $g$ is merely lsc but $\nabla f$ is locally Lipschitz continuous.
		From $x^k \to_K x^\star$ and lsc of $g$, we find
		\begin{equation*}
			g (x^\star) \leq \liminf_{k \in K} g (x^k) \leq \limsup_{k \in K} g (x^k) .
		\end{equation*}
		Therefore, it suffices to show that $\limsup_{k \in K} g (x^k) \leq g(x^\star)$.
		Since $x^k$ solves the subproblem at \cref{state:NMPG:barx} with stepsize $\gamma_k$, we obtain
		\begin{equation*}
			\innprod{\nabla f (x^{k-1})}{x^k - x^\star} + \frac{1}{2 \gamma_k} \| x^k - x^{k-1} \|^2 +  g (x^k )
			\leq
			\frac{1}{2 \gamma_k} \| x^\star - x^{k-1} \|^2 + g (x^\star)
		\end{equation*}
		for each $k \in K$.
		We now take the upper limit for $k \to_K \infty$ on both sides.
		By \cref{cor:gammaLocalLipschitz}, the stepsizes $\seq{\gamma_k}{k\in K}$ remain bounded away from zero, hence $\frac{1}{\gamma_k} \| x^\star - x^{k-1} \|^2 \to_K 0$.
		Using the continuity of $\nabla f$, the convergences $x^k - x^{k-1} \to_K 0$ as well as $\frac{1}{\gamma_k} \| x^k - x^{k-1} \|^2 \to_K 0$ (see \cref{thm:asymp}\ref{thm:asymp:summable}-\ref{thm:asymp:diff_x}), we obtain $\limsup_{k \in K} g(x^k) \leq g(x^\star)$.
	\end{enumerate}
	Altogether, we therefore get $x^k \to_K x^\star$ with $\varphi(x^k) \to_K \varphi(x^\star) = \varphi_\star$, completing the proof.
\end{proof}
We accompany our asymptotic characterization of \cref{alg:NMPG} with the finite termination property.
\begin{mybox}
	\begin{theorem}\label{thm:finiteTermination}
		Suppose that \crefrange{ass:f}{ass:phi} hold.
		Let any $x^0 \in \dom g$ and $\varepsilon > 0$ be provided.
		If the sequence $\seq{x^k}{k\in \N}$ generated by \cref{alg:NMPG} admits an accumulation point, then the algorithm returns an $\varepsilon$-M-stationary point of \cref{eq:P} after finitely many iterations.
	\end{theorem}
\end{mybox}
\begin{proof}
	Let $x^\star \in \XX$ be an accumulation point of $\seq{ x^k }{k\in\N}$ and let $\seq{ x^k }{k\in K}$ be a subsequence such that $x^k \to_K x^\star$.
	Then, we have that $x^{k-1} \to_K x^\star$ by \cref{thm:asymp}\ref{thm:asymp:diff_x} and $\frac{1}{\gamma_k} (x^k - x^{k-1}) \to_K 0$ by \cref{thm:asymp:fpr}.
	Thus, owing to the continuity of $\nabla f$ by \cref{ass:f}, the stopping condition at \cref{state:NMPG:termination} is satisfied for some $k \in K$ large enough.
	Finally, $\varepsilon$-M-stationarity of the returned point $x^k$ follows from the bounds in \cref{eq:boundDistanceSubdifferential}.
	\qedhere
\end{proof}
	
	\subsection{Convergence Rates}\label{sec:ConvergenceRate}
	We now turn to analyzing the rate of asymptotic regularity for \cref{alg:NMPG} and, considering the iterates it generates, we recover the classical (worst-case) rate of (at least) $\bigO(1/\sqrt{k})$ for nonconvex problems \cite[Thm 10.15(c)]{beck2017first}.
Noticing that $\gamma_k \leq \gamma_{\max}$ for all $k$ and $\Phi_0 = \varphi(x^0) < \infty$, the result readily follows from the inequalities in \cref{eq:NMPG:asymp:summable}.
A detailed proof is included in the Additional Material (p.\ \pageref{proof:thm:asymp:localRateIterates}).
\begin{mybox}
	\begin{corollary}\label{thm:asymp:localRateIterates}%
		Suppose that \crefrange{ass:f}{ass:phi} are satisfied and consider a sequence of iterates $\seq{ x^k }{k\in\N}$ generated by \cref{alg:NMPG}.
		Then, for every $k \geq 1$ we have
		\begin{equation*}
			\min_{1 \leq j \leq k} \| x^j - x^{j-1} \|
			\leq
			\frac{1}{\sqrt{k}}
			\left( \frac{2 \gamma_{\max} [\varphi(x^0) - \varphi_\star]}{p_{\min} (1 - \alpha)} \right)^{1/2}
			<
			\infty
		\end{equation*}
		where $\varphi_\star$ is that of \cref{thm:asymp}\ref{thm:asymp:cost}.
	\end{corollary}
\end{mybox}
We are also interested in characterizing the convergence rate of some stationarity measure for \cref{eq:P}.
In the nonconvex setting, we shall monitor the quantity $r_k \coloneqq \frac{1}{\gamma_k} \| x^k - x^{k-1} \|$ as a measure of progress within \cref{alg:NMPG}.
Let us denote the proximal gradient mapping by $\ffunc{\Tprox_\gamma}{\XX}{\XX}$, $\Tprox_{\gamma}( x ) \coloneqq \prox_{\gamma g}(x - \gamma \nabla f(x))$.
Since $x^k \in \Tprox_{\gamma_k}( x^{k-1} )$ by \cref{state:NMPG:barx}, we have that $\dist(0,\Rprox_{\gamma_k}(x^{k-1})) \leq r_k$, where $\ffunc{\Rprox_\gamma}{\XX}{\XX}$, $\Rprox_{\gamma}(x) \coloneqq \frac{1}{\gamma} [x - \Tprox_\gamma(x)]$, denotes the fixed-point residual associated to \cref{eq:P} and provides a measure of stationarity thereof \cite{beck2017first,stella2017simple,themelis2018forward}.

In order to obtain global complexity results, it seems indispensable the assumption of local Lipschitz continuity of $\nabla f$.
In fact, revealing that the stepsizes $\gamma_k$ remain bounded away from zero (possibly only subsequentially or under some compactness assumptions), \cref{cor:gammaLocalLipschitz} plays a central role in obtaining suitable bounds for inferring rate-of-convergence results without additional assumptions.

We begin by stating in \cref{thm:localLipschitzSubsequenceSublinearRate} that the global convergence rate of a stationarity measure for \cref{eq:P} is (at least) sublinear for convergent \emph{subsequences}.
\begin{mybox}
	\begin{theorem}\label{thm:localLipschitzSubsequenceSublinearRate}
		Suppose \crefrange{ass:f}{ass:phi} are satisfied and $\nabla f$ is locally Lipschitz continuous.
		Consider a sequence $\seq{x^k}{k\in\N}$ generated by \cref{alg:NMPG}.
		Let $x^\star$ be an accumulation point of $\seq{ x^k }{k\in\N}$ and $\seq{x^k}{k \in K}$ a subsequence such that $x^k \to_K x^\star$.
		Then for every $k \geq 1$ we have
		\begin{equation*}
			\min_{j \in K_k} \frac{\| x^j - x^{j-1} \|}{\gamma_j}
			\leq
			\frac{1}{\sqrt{\cardinality{K_k}}}
			\left( \frac{2 [\varphi(x^0) - \varphi_\star]}{\gamma_\star p_{\min} (1 - \alpha)} \right)^{1/2}
		\end{equation*}
		where $\varphi_\star$ is that of \cref{thm:asymp}\ref{thm:asymp:cost}, $\gamma_\star > 0$ is a lower bound on $\seq{\gamma_k}{k \in K}$, and $\cardinality{K_k}$ denotes the cardinality of the set of indices $K_k \coloneqq \{j \in K \,\vert\, 1 \leq j \leq k \}$.
	\end{theorem}
\end{mybox}
\begin{proof}
	Since $\nabla f$ is locally Lipschitz continuous and $x^k \to_K x^\star$, the subsequence $\seq{\gamma_k}{k\in K}$ is bounded away from zero by some $\gamma_\star > 0$, owing to \cref{cor:gammaLocalLipschitz}.
	For all $k \in \N$, let $r_k \coloneqq \| x^k - x^{k-1} \| / \gamma_k$.
	Extracting from \cref{eq:NMPG:asymp:summable}, using that $\gamma_k \geq \gamma_\star > 0$ for all $k \in K$, and rearranging, we obtain
	\begin{equation*}
		\frac{2 [\Phi_0 - \varphi_\star]}{p_{\min} (1 - \alpha)}
		\geq
		\sum_{j=1}^k \gamma_j r_j^2
		\geq
		\sum_{j \in K_k} \gamma_j r_j^2
		\geq
		\gamma_\star \sum_{j \in K_k} r_j^2
		\geq
		\gamma_\star \cardinality{K_k} \min_{j \in K_k} r_j^2 ,
	\end{equation*}
	proving the result, since $\Phi_0 = \varphi(x^0) < \infty$.
\end{proof}
Now, if the iterates $\seq{x^k}{k\in\N}$ remain bounded, local Lipschitzness and compactness yield a uniform lower bound on the stepsizes $\seq{\gamma_k}{k\in\N}$.
Then, a global sublinear rate of convergence for the stationarity measure readily follows from that of the iterates, given in \cref{thm:asymp:localRateIterates}.
However, boundedness of the iterates $\seq{x^k}{k\in\N}$ can be a strong assumption.
In the following \cref{thm:localLipschitzAccumulationCompactBoundedGamma}, we demonstrate that boundedness of the set of accumulation points is sufficient for the stepsizes to remain uniformly away from zero.
\begin{mybox}
	\begin{lemma}\label{thm:localLipschitzAccumulationCompactBoundedGamma}
		Suppose \crefrange{ass:f}{ass:phi} are satisfied and $\nabla f$ is locally Lipschitz continuous.
		Consider a sequence $\seq{x^k}{k\in\N}$ generated by \cref{alg:NMPG}.
		If the set of accumulation points of $\seq{x^k}{k\in\N}$ is bounded, then there exists $\overline{\gamma} > 0$ such that $\gamma_k \geq \overline{\gamma}$ for all $k\in\N$.
	\end{lemma}
\end{mybox}
\begin{proof}
	Let $\omega \subset \XX$ denote the set of accumulation points of $\seq{x^k}{k\in\N}$.
	Owing to \cref{thm:asymp}\ref{thm:asymp:diff_x}, if $\omega$ is bounded it is also compact.
	Now suppose, \emph{ad absurdum}, that $\gamma_k \to 0$.
	Then, by compactness of $\omega$, there exists a subsequence $\seq{x^k}{k \in K}$ such that $x^k \to_K x^\star$ for some $x^\star \in \omega$ and at the same time $\gamma_k \to_K 0$.
	On the contrary, due to the locally Lipschitz continuity assumption, \cref{cor:gammaLocalLipschitz} prevents $\seq{\gamma_k}{k \in K}$ from vanishing, thus proving the claim by contradiction.
\end{proof}
Finally, patterning \cref{thm:localLipschitzSubsequenceSublinearRate}, the following result gives the corresponding global convergence rate with respect to the whole sequence.
A detailed proof is included in the Additional Material (p.\ \pageref{proof:thm:localLipschitzAccumulationCompactSublinearRate}).
\begin{mybox}
	\begin{corollary}\label{thm:localLipschitzAccumulationCompactSublinearRate}
		Suppose \crefrange{ass:f}{ass:phi} are satisfied and $\nabla f$ is locally Lipschitz continuous.
		Consider a sequence $\seq{x^k}{k\in\N}$ generated by \cref{alg:NMPG}.
		If the set of accumulation points of $\seq{x^k}{k\in\N}$ is compact, then for every $k \geq 1$ we have
		\begin{equation*}
			\min_{1 \leq j \leq k} \frac{\| x^j - x^{j-1} \|}{\gamma_j}
			\leq
			\frac{1}{\sqrt{k}} \left( \frac{2 [\varphi(x^0) - \varphi_\star]}{\gamma_\star p_{\min} (1 - \alpha)} \right)^{1/2}
		\end{equation*}
		where $\varphi_\star$ is that of \cref{thm:asymp}\ref{thm:asymp:cost} and $\gamma_\star > 0$ is a lower bound for $\seq{\gamma_k}{k \in \N}$.
	\end{corollary}
\end{mybox}
Similar results can be obtained if the entire sequence $\seq{x^k}{k\in\N}$ converges to some accumulation point $x^\ast$.
A situation leading to this circumstance is that of $\varphi$ being locally strongly convex in a neighborhood of $x^\ast$.
In this case, the convergence of $\seq{x^k}{k\in\N}$ to $x^\ast$ directly follows from \cref{thm:asymp}\ref{thm:asymp:diff_x} and \cref{thm:NMPG:subseqStationarity}, owing to \cite[Lemma 4.10]{more1983computing}.
Another occurrence is that of $\varphi$ satisfying some local error bound conditions \cite[Thm 2]{tseng2009coordinate}, closely related to the Kurdyka-{\L}ojasiewicz property \cite{li2018calculus}.

\section{Spectral Stepsize and Numerical Example}\label{sec:NumericalExample}
The general scheme we investigated comprises an arbitrary (yet bounded) choice of the proximal stepsize, enables nonmonotone behavior, and provides convergence guarantees under weak assumptions.
These features lay the foundation for \emph{spectral proximal gradient methods}, which can be constructed by considering the Barzilai-Borwein spectral estimate as a principled way for selecting the proximal stepsize \cite{barzilai1988two}.
In the spirit of quasi-Newton methods, a two-point approximation yields the spectral stepsize obtained in \cite[Eq. 6]{barzilai1988two}:
\begin{equation}\label{eq:spectral_stepsize}
\gamma_k^{\text{BB}} \coloneqq \frac{\innprod{\Delta x^k}{\Delta x^k}}{\innprod{\Delta x^k}{\Delta g^k}}
\qquad\text{with}~
\Delta x^k \coloneqq x^k - x^{k-1}
\quad\text{and}\quad
\Delta g^k \coloneqq \nabla f(x^k) - \nabla f(x^{k-1}) .
\end{equation}
Associated to a local (isotropic, quadratic) surrogate of the smooth term $f$, the spectral choice is the ``essential feature that puts efficiency in the projected gradient methodology'' \cite[\S 3]{birgin2000nonmonotone}.
When coupled with a linesearch for globalization, nonmonotonicity proved to be key in unleashing the performance of projected gradient methods with spectral approximations \cite{raydan1997barzilai,birgin2014spectral}.
We hope this behavior carries over to proximal methods.

We now compare the practical performance of several variants of \cref{alg:NMPG} on a numerical example.
Two stepsize selection strategies are considered: ``plain'' simply inherits the previous value $\gamma_{k-1}$, whereas ``spectral'' uses the estimate $\gamma_{k-1}^{\text{BB}}$ in \eqref{eq:spectral_stepsize}; both choices are initialized with $\gamma_0 \coloneqq 1$ and projected onto $[\gamma_{\min}, \gamma_{\max}]$ to comply with \cref{state:NMPG:init}.
These strategies are combined with three linesearch flavors: ``monotone'' ($p_k \coloneqq 1$ for all $k\in \N$), nonmonotone with ``average'' merit ($p_k \coloneqq 0.2$ for all $k\in\N$), and nonmonotone with ``max'' merit (memory $M = 5$) as in \eqref{eq:nm_flavours}.
Other parameters of \cref{alg:NMPG} are set to $\alpha_k \coloneqq 0.999$, $\beta_k \coloneqq 0.5$ for all $k\in\N$, $\gamma_{\min} \coloneqq 10^{-12}$, $\gamma_{\max} \coloneqq 10^{12}$, and tolerance $\varepsilon \coloneqq 10^{-6}$.

\subsection*{Dictionary Learning}

The challenge of dictionary learning is that of finding a collection of vectors that can sparsely yet accurately represent data signals.
Given $m$ signals $y_1,\ldots, y_m \in \R^n$, we wish to find $\ell$ dictionary atoms $d_1, \ldots, d_\ell \in \R^n$ such that each signal $y_j$ can be approximated by a linear combination of these atoms via a sparse vector of coefficients $c_j \in \R^\ell$.
Stacking the problem data in $Y \in \R^{n \times m}$, the unknown dictionary atoms in $D \in \R^{n \times \ell}$, and the unknown coefficients in $C \in \R^{\ell \times m}$, the problem can be expressed as
\begin{align}\label{eq:dictionarylearning}
	\minimize_{D,C} \quad& \frac{1}{2} \| Y - D C \|_F^2 + \lambda \| C \|_0 &
	\stt \quad& \| d_j \|_2 = 1 \quad j=1,\ldots,\ell ,
\end{align}
where $\| \cdot \|_F$ and $\| \cdot \|_2$ denote the Frobenius and Euclidean norm, respectively.
Parameter $\lambda \geq 0$ is introduced to tune the sparsity-promoting regularization $\| C \|_0$, which counts the number of nonzero entries in $C$.
Without loss of generality, the dictionary atoms $d_j$ are constrained to have unit norm; by scaling the coefficients in $C$, the objective remains unchanged.
Relative to \eqref{eq:P}, function $g$ encapsulates both the norm constraints on $D$ and the cardinality cost on $C$, while $f$ gives the quadratic loss.

We tested our algorithm on $100$ problem instances with $n = 10$, $\ell = 20$ and $m = 30$, generated as follows.
First, a dictionary $D \in \R^{n \times \ell}$ was constructed with entries sampled from a standard normal distribution, and each column was scaled to have unit norm.
Then, a matrix $C \in \R^{\ell \times m}$ was randomly generated with $N = 3$ normally distributed nonzero coefficients per column. 
Finally, we set $Y \coloneqq D C$, considered $\lambda \coloneqq 10^{-2}$, and initialized the algorithm with $D^0 \in \R^{n \times \ell}$ and $C^0 \in \R^{\ell \times m}$ filled with normally distributed entries.

\begin{figure}
	\centering
	\definecolor{mycolor1}{rgb}{0.89412,0.10196,0.10980}%
\definecolor{mycolor2}{rgb}{0.30196,0.68627,0.29020}%
\definecolor{mycolor3}{rgb}{0.21569,0.49412,0.72157}%

\begin{tikzpicture}

\begin{axis}[%
width=0.43\columnwidth,
height=0.266\columnwidth,
at={(0,0)},
scale only axis,
axis x line=bottom,
axis y line=left,
xmode=log,
xmin=300,
xmax=200000,
xlabel style={font=\color{white!15!black}},
xlabel={Number of proximal evaluations},
ymin=-0.05,
ymax=1.05,
ylabel style={font=\color{white!15!black}},
ylabel={Fraction of problems solved},
]

\addplot [color=mycolor3, line width=1.4pt, forget plot]
table {TeX/Tikz/Data/dictionary_learning_lambda001_prox/spectral_monotone.dat};

\addplot [color=mycolor2, line width=1.4pt, forget plot]
table {TeX/Tikz/Data/dictionary_learning_lambda001_prox/spectral_max.dat};

\addplot [color=mycolor1, line width=1.4pt, forget plot]
table {TeX/Tikz/Data/dictionary_learning_lambda001_prox/spectral_average.dat};

\addplot [color=mycolor3, dotted, line width=1.4pt, forget plot]
table {TeX/Tikz/Data/dictionary_learning_lambda001_prox/plain_monotone.dat};

\addplot [color=mycolor2, dotted, line width=1.4pt, forget plot]
table {TeX/Tikz/Data/dictionary_learning_lambda001_prox/plain_max.dat};

\addplot [color=mycolor1, dotted, line width=1.4pt, forget plot]
table {TeX/Tikz/Data/dictionary_learning_lambda001_prox/plain_average.dat};

\end{axis}

\begin{axis}[%
width=0.43\columnwidth,
height=0.266\columnwidth,
at={(0.47\columnwidth,0)},
scale only axis,
axis x line=bottom,
axis y line=left,
xmin=4.2,
xmax=6.1,
xlabel style={font=\color{white!15!black}},
xlabel={Objective},
ymin=-0.05,
ymax=1.05,
legend style={at={(0.97,0.03)}, anchor=south east, legend cell align=left, align=left, draw=white}
]

\addplot [color=mycolor3, line width=1.4pt]
table {TeX/Tikz/Data/dictionary_learning_lambda001_obj/spectral_monotone.dat};
\addlegendentry{monotone}

\addplot [color=mycolor2, line width=1.4pt]
table {TeX/Tikz/Data/dictionary_learning_lambda001_obj/spectral_max.dat};
\addlegendentry{max}

\addplot [color=mycolor1, line width=1.4pt]
table {TeX/Tikz/Data/dictionary_learning_lambda001_obj/spectral_average.dat};
\addlegendentry{average}

\addplot [color=mycolor3, dotted, line width=1.4pt, forget plot]
table {TeX/Tikz/Data/dictionary_learning_lambda001_obj/plain_monotone.dat};

\addplot [color=mycolor2, dotted, line width=1.4pt, forget plot]
table {TeX/Tikz/Data/dictionary_learning_lambda001_obj/plain_max.dat};

\addplot [color=mycolor1, dotted, line width=1.4pt, forget plot]
table {TeX/Tikz/Data/dictionary_learning_lambda001_obj/plain_average.dat};

\end{axis}

\end{tikzpicture}%
	\caption{%
		Comparing variants of \cref{alg:NMPG} on dictionary learning \eqref{eq:dictionarylearning}: combinations of plain (dotted) and spectral (solid) stepsizes with monotone and nonmonotone linesearch strategies.%
	}%
	\label{fig:dictionarylearning}
\end{figure}
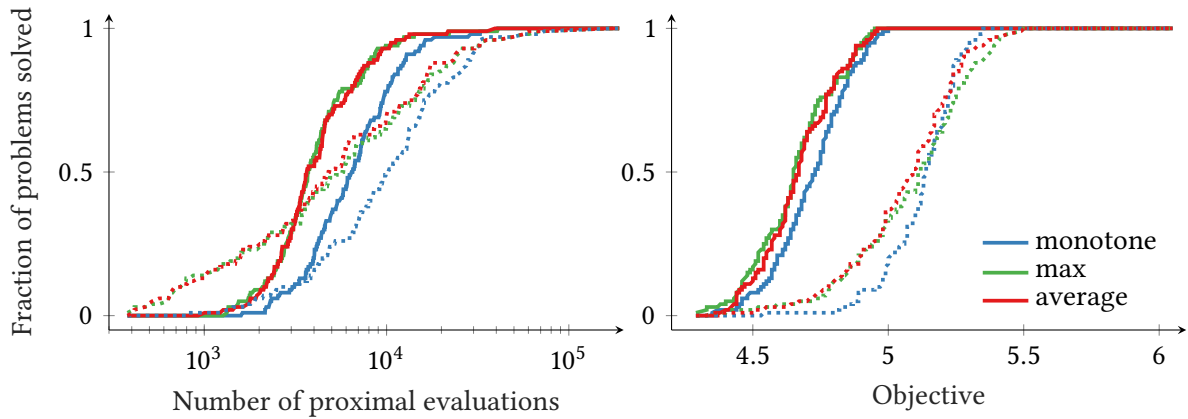

Numerical results are depicted in \cref{fig:dictionarylearning}, using profiles that represent the percentage of problems solved as a function of the available computational budget, considering either the evaluations of $\prox_g$ or the final value of the objective $\varphi$.
Note that all solver variants find a solution, in the sense of \eqref{eq:boundDistanceSubdifferential}, and terminate.
Possibly inducing more erratic behaviors, spectral stepsize selection tends to reach lower objective values, often with fewer evaluations compared to the plain counterpart.
Conversely, for both plain and spectral stepsizes, nonmonotone strategies appear to perform better than monotone variants, especially in terms of evaluations (function, proximal, and iterations).
Although the max and averaging nonmonotonicity flavors exhibit comparable capabilities, guarantees under weak assumptions are currently available for the latter only.

\section{Conclusions}\label{sec:Conclusions}
We studied an adaptive nonmonotone proximal gradient method for fully nonconvex composite optimization and its convergence properties under weak assumptions.
Our findings indicate that, by attaining sufficient decrease at every iteration, an averaging merit function allows to establish convergence guarantees in the nonmonotone regime, matching those in the monotone case and thus bridging the gap emerged in the literature \cite{kanzow2022convergence}.
Furthermore, by exploiting these characterizations, we derived global rate-of-convergence results and bounds, for both iterates and a stationarity measure.
Numerical comparisons demonstrated practical benefits of embracing nonmonotone linesearch strategies, whose averaging variant can rely on theoretical guarantees.

Future research work shall study the trade-offs and further compare (the performance of) spectral proximal gradient methods against Newton-type and variable metric approaches \cite{lee2014proximal}, which are based on more detailed (but potentially more intricate) models of the local landscape.
Moreover, it remains an open question whether PANOC-type methods \cite{stella2017simple,demarchi2022proximal} maintain or lose their convergence properties when the smooth term $f$ is merely continuously differentiable.

\section*{Acknowledgments}
\TheAcknowledgments

\phantomsection
\addcontentsline{toc}{section}{Additional Material}%
\section*{Additional Material}
\begin{proof}[Proof of \cref{cor:gammaLocalLipschitz}]\label{proof:cor:gammaLocalLipschitz}
	We may argue as in the proof of \cref{thm:asymp:fpr}.
	By contradiction, let us assume, without loss of generality, that $\gamma_k \to_K 0$.
	For each $k \in K$, define $\hat{\gamma}_k$ and $\hat{x}^k$ as in that proof, and let $L_\star > 0$ denote the local Lipschitz constant of $\nabla f$ in a neighborhood of $x^\star$.
	Recall that $x^k \to_K x^\star$ and, from the proof of \cref{thm:asymp:fpr}, that $\hat{x}^k \to_K x^\star$.
	Exploiting \cref{Prop:gammaxdiff_NM-4}, we therefore obtain
	\begin{equation*}
		\frac{\alpha}{2 \hat{\gamma}_k} \| \hat{x}^k - x^{k-1} \|
		{}\leq{}
		\| \nabla f(x^{k-1}) - \nabla f(\xi^k) \|
		{}\leq{}
		L_\star \| x^{k-1} - \xi^k \|
		{}\leq{}
		L_\star \| \hat{x}^k - x^{k-1} \|
	\end{equation*}
	for all $k \in K$ sufficiently large, using the fact that $\xi^k$ is on the line segment between $x^{k-1}$ and $\hat{x}^k$.
	Since $\hat{x}^k \neq x^{k-1}$ by \cref{Prop:gammaxdiff_NM-1}, this gives a contradiction as $\hat{\gamma}_k \to_K 0$. 
	Hence, $\seq{ \gamma_k }{k\in K}$ stays bounded away from zero.
\end{proof}

\begin{proof}[Proof of \cref{thm:asymp:localRateIterates}]\label{proof:thm:asymp:localRateIterates}
	Let us denote $\Delta_k \coloneqq \| x^k - x^{k-1} \|$ for all $k \geq 1$.
	Owing to $\gamma_k \leq \gamma_{\max}$ for all $k$ and $\Phi_0 = \varphi(x^0) < \infty$, from \cref{eq:NMPG:asymp:summable} we obtain that
	\begin{equation*}
		\frac{2 [\varphi(x^0) - \varphi_\star]}{p_{\min} (1 - \alpha)}
		\geq
		\sum_{j=1}^k \frac{1}{\gamma_j} \Delta_j^2
		\geq
		\frac{1}{\gamma_{\max}} \sum_{j=1}^k \Delta_j^2
		\geq
		\frac{k}{\gamma_{\max}} \min_{1 \leq j \leq k} \Delta_j^2 .
	\end{equation*}
	Rearranging yields the claimed assertion.
\end{proof}

\begin{proof}[Proof of \cref{thm:localLipschitzAccumulationCompactSublinearRate}]\label{proof:thm:localLipschitzAccumulationCompactSublinearRate}
	Since $\nabla f$ is locally Lipschitz continuous, it follows from \cref{thm:localLipschitzAccumulationCompactBoundedGamma} that $\seq{\gamma_k}{k\in\N}$ remains bounded away from zero if the set of accumulation points of $\seq{x^k}{k\in\N}$ is compact.
	Rearranging from \cref{eq:NMPG:asymp:summable}, using that $\gamma_k \geq \gamma_\star > 0$ for all $k \in \N$, and rearranging, we obtain
	\begin{equation*}
		\frac{2 [\Phi_0 - \varphi_\star]}{p_{\min} (1 - \alpha)}
		\geq
		\sum_{j=1}^k \frac{1}{\gamma_j} \| x^j - x^{j-1} \|^2
		\geq
		\gamma_\star \sum_{j = 1}^k \frac{\| x^j - x^{j-1} \|^2}{\gamma_j^2}
		\geq
		k \gamma_\star \min_{1 \leq j \leq k} \frac{\| x^j - x^{j-1} \|^2}{\gamma_j^2} ,
	\end{equation*}
	proving the result.
\end{proof}

\bibliographystyle{jnsao}
\bibliography{spectralproxgrad}

\begin{thebibliography}{10}

\bibitem{armijo1966minimization}
L{.\nobreak\kern 0.33333em}Armijo, Minimization of functions having {L}ipschitz
  continuous first partial derivatives, \emph{Pacific Journal of Mathematics}
  16 (1966),  1--3.

\bibitem{barzilai1988two}
J{.\nobreak\kern 0.33333em}Barzilai and J.\,M{.\nobreak\kern 0.33333em}Borwein,
  Two-point step size gradient methods, \emph{IMA Journal of Numerical
  Analysis} 8 (1988),  141--148.

\bibitem{bauschke2017convex}
H.\,H{.\nobreak\kern 0.33333em}Bauschke and P.\,L{.\nobreak\kern
  0.33333em}Combettes, \emph{Convex Analysis and Monotone Operator Theory in
  {H}ilbert Spaces}, Springer International Publishing, 2017,
  \href{https://dx.doi.org/10.1007/978-3-319-48311-5}{\nolinkurl{doi:10.1007/978-3-319-48311-5}}.

\bibitem{beck2017first}
A{.\nobreak\kern 0.33333em}Beck, \emph{First-Order Methods in Optimization},
  Society for Industrial and Applied Mathematics, Philadelphia, PA, 2017,
  \href{https://dx.doi.org/10.1137/1.9781611974997}{\nolinkurl{doi:10.1137/1.9781611974997}}.

\bibitem{birgin2014practical}
E.\,G{.\nobreak\kern 0.33333em}Birgin and J.\,M{.\nobreak\kern
  0.33333em}Mart{\'i}nez, \emph{Practical Augmented {L}agrangian Methods for
  Constrained Optimization}, Society for Industrial and Applied Mathematics,
  Philadelphia, PA, 2014.

\bibitem{birgin2000nonmonotone}
E.\,G{.\nobreak\kern 0.33333em}Birgin, J.\,M{.\nobreak\kern
  0.33333em}Mart{\'i}nez, and M{.\nobreak\kern 0.33333em}Raydan, Nonmonotone
  spectral projected gradient methods on convex sets, \emph{SIAM Journal on
  Optimization} 10 (2000),  1196--1211,
  \href{https://dx.doi.org/10.1137/S1052623497330963}{\nolinkurl{doi:10.1137/s1052623497330963}}.

\bibitem{birgin2014spectral}
E.\,G{.\nobreak\kern 0.33333em}Birgin, J.\,M{.\nobreak\kern
  0.33333em}Mart{\'i}nez, and M{.\nobreak\kern 0.33333em}Raydan, Spectral
  projected gradient methods: Review and perspectives, \emph{Journal of
  Statistical Software} 60 (2014),  1--21,
  \href{https://dx.doi.org/10.18637/jss.v060.i03}{\nolinkurl{doi:10.18637/jss.v060.i03}}.

\bibitem{cauchy1847methode}
A.\,L{.\nobreak\kern 0.33333em}Cauchy, M\'ethode g\'en\'erale pour la
  r\'esolution des syst\`ems d'\'equations simultan\'ees, \emph{Comp. Rend.
  Sci. Paris} 25 (1847),  536--538.

\bibitem{combettes2011proximal}
P.\,L{.\nobreak\kern 0.33333em}Combettes and J.\,C{.\nobreak\kern
  0.33333em}Pesquet, Proximal splitting methods in signal processing, in
  \emph{Fixed-Point Algorithms for Inverse Problems in Science and
  Engineering}, Springer, New York, 2011,  185--212.

\bibitem{demarchi2023implicit}
A{.\nobreak\kern 0.33333em}De~Marchi, Implicit augmented {L}agrangian and
  generalized optimization, 2023,
  \href{https://arxiv.org/abs/2302.00363}{\nolinkurl{arXiv:2302.00363}}.

\bibitem{demarchi2023constrained}
A{.\nobreak\kern 0.33333em}De~Marchi, X{.\nobreak\kern 0.33333em}Jia,
  C{.\nobreak\kern 0.33333em}Kanzow, and P{.\nobreak\kern 0.33333em}Mehlitz,
  Constrained composite optimization and augmented {L}agrangian methods,
  \emph{Mathematical Programming}  (2023),
  \href{https://dx.doi.org/10.1007/s10107-022-01922-4}{\nolinkurl{doi:10.1007/s10107-022-01922-4}}.

\bibitem{demarchi2022interior}
A{.\nobreak\kern 0.33333em}De~Marchi and A{.\nobreak\kern 0.33333em}Themelis,
  An interior proximal gradient method for nonconvex optimization, 2022,
  \href{https://arxiv.org/abs/2208.00799}{\nolinkurl{arXiv:2208.00799}}.

\bibitem{demarchi2022proximal}
A{.\nobreak\kern 0.33333em}De~Marchi and A{.\nobreak\kern 0.33333em}Themelis,
  Proximal gradient algorithms under local {L}ipschitz gradient continuity,
  \emph{Journal of Optimization Theory and Applications} 194 (2022),  771--794,
  \href{https://dx.doi.org/10.1007/s10957-022-02048-5}{\nolinkurl{doi:10.1007/s10957-022-02048-5}}.

\bibitem{fukushima1981generalized}
M{.\nobreak\kern 0.33333em}Fukushima and H{.\nobreak\kern 0.33333em}Mine, A
  generalized proximal point algorithm for certain non-convex minimization
  problems, \emph{International Journal of Systems Science} 12 (1981),
  989--1000,
  \href{https://dx.doi.org/10.1080/00207728108963798}{\nolinkurl{doi:10.1080/00207728108963798}}.

\bibitem{grippo1986nonmonotone}
L{.\nobreak\kern 0.33333em}Grippo, F{.\nobreak\kern 0.33333em}Lampariello, and
  S{.\nobreak\kern 0.33333em}Lucidi, A nonmonotone line search technique for
  {N}ewton's method, \emph{SIAM Journal on Numerical Analysis} 23 (1986),
  707--716,
  \href{https://dx.doi.org/10.1137/0723046}{\nolinkurl{doi:10.1137/0723046}}.

\bibitem{jia2022augmented}
X{.\nobreak\kern 0.33333em}Jia, C{.\nobreak\kern 0.33333em}Kanzow,
  P{.\nobreak\kern 0.33333em}Mehlitz, and G{.\nobreak\kern 0.33333em}Wachsmuth,
  An augmented {L}agrangian method for optimization problems with structured
  geometric constraints, \emph{Mathematical Programming}  (2022),
  \href{https://dx.doi.org/10.1007/s10107-022-01870-z}{\nolinkurl{doi:10.1007/s10107-022-01870-z}}.

\bibitem{kanzow2022convergence}
C{.\nobreak\kern 0.33333em}Kanzow and P{.\nobreak\kern 0.33333em}Mehlitz,
  Convergence properties of monotone and nonmonotone proximal gradient methods
  revisited, \emph{Journal of Optimization Theory and Applications}  (2022),
  \href{https://dx.doi.org/10.1007/s10957-022-02101-3}{\nolinkurl{doi:10.1007/s10957-022-02101-3}}.

\bibitem{lee2014proximal}
J.\,D{.\nobreak\kern 0.33333em}Lee, Y{.\nobreak\kern 0.33333em}Sun, and
  M.\,A{.\nobreak\kern 0.33333em}Saunders, Proximal {N}ewton-type methods for
  minimizing composite functions, \emph{SIAM Journal on Optimization} 24
  (2014),  1420--1443,
  \href{https://dx.doi.org/10.1137/130921428}{\nolinkurl{doi:10.1137/130921428}}.

\bibitem{li2018calculus}
G{.\nobreak\kern 0.33333em}Li and T.\,K{.\nobreak\kern 0.33333em}Pong, Calculus
  of the exponent of Kurdyka-{\L}ojasiewicz inequality and its applications to
  linear convergence of first-order methods, \emph{Foundations of Computational
  Mathematics} 18 (2018),  1199--1232,
  \href{https://dx.doi.org/10.1007/s10208-017-9366-8}{\nolinkurl{doi:10.1007/s10208-017-9366-8}}.

\bibitem{li2015accelerated}
H{.\nobreak\kern 0.33333em}Li and Z{.\nobreak\kern 0.33333em}Lin, Accelerated
  proximal gradient methods for nonconvex programming, in \emph{Advances in
  Neural Information Processing Systems}, C{.\nobreak\kern 0.33333em}Cortes,
  N{.\nobreak\kern 0.33333em}Lawrence, D{.\nobreak\kern 0.33333em}Lee,
  M{.\nobreak\kern 0.33333em}Sugiyama, and R{.\nobreak\kern 0.33333em}Garnett
  (eds.), volume~28, Curran Associates, Inc., 2015.

\bibitem{mordukhovich2018variational}
B.\,S{.\nobreak\kern 0.33333em}Mordukhovich, \emph{Variational Analysis and
  Applications}, Springer, 2018,
  \href{https://dx.doi.org/10.1007/978-3-319-92775-6}{\nolinkurl{doi:10.1007/978-3-319-92775-6}}.

\bibitem{more1983computing}
J.\,J{.\nobreak\kern 0.33333em}Mor{\'e} and D.\,C{.\nobreak\kern
  0.33333em}Sorensen, Computing a trust region step, \emph{SIAM Journal on
  Scientific and Statistical Computing} 4 (1983),  553--572,
  \href{https://dx.doi.org/10.1137/0904038}{\nolinkurl{doi:10.1137/0904038}}.

\bibitem{parikh2014proximal}
N{.\nobreak\kern 0.33333em}Parikh and S{.\nobreak\kern 0.33333em}Boyd, Proximal
  algorithms, \emph{Foundations and Trends in Optimization} 1 (2014),
  127--239,
  \href{https://dx.doi.org/10.1561/2400000003}{\nolinkurl{doi:10.1561/2400000003}}.

\bibitem{raydan1997barzilai}
M{.\nobreak\kern 0.33333em}Raydan, The {B}arzilai and {B}orwein gradient method
  for the large scale unconstrained minimization problem, \emph{SIAM Journal on
  Optimization} 7 (1997),  26--33,
  \href{https://dx.doi.org/10.1137/S1052623494266365}{\nolinkurl{doi:10.1137/s1052623494266365}}.

\bibitem{rockafellar1998variational}
R.\,T{.\nobreak\kern 0.33333em}Rockafellar and R.\,J{.\nobreak\kern
  0.33333em}Wets, \emph{Variational Analysis}, volume 317, Springer, 1998.

\bibitem{sopasakis2020open}
P{.\nobreak\kern 0.33333em}Sopasakis, E{.\nobreak\kern 0.33333em}Fresk, and
  P{.\nobreak\kern 0.33333em}Patrinos, {OpEn}: Code generation for embedded
  nonconvex optimization, \emph{IFAC-PapersOnLine} 53 (2020),  6548--6554,
  \href{https://dx.doi.org/10.1016/j.ifacol.2020.12.071}{\nolinkurl{doi:10.1016/j.ifacol.2020.12.071}}.

\bibitem{stella2017simple}
L{.\nobreak\kern 0.33333em}Stella, A{.\nobreak\kern 0.33333em}Themelis,
  P{.\nobreak\kern 0.33333em}Sopasakis, and P{.\nobreak\kern
  0.33333em}Patrinos, A simple and efficient algorithm for nonlinear model
  predictive control, in \emph{2017 IEEE 56th Annual Conference on Decision and
  Control (CDC)}, IEEE, 2017,  1939--1944,
  \href{https://dx.doi.org/10.1109/CDC.2017.8263933}{\nolinkurl{doi:10.1109/cdc.2017.8263933}}.

\bibitem{themelis2018forward}
A{.\nobreak\kern 0.33333em}Themelis, L{.\nobreak\kern 0.33333em}Stella, and
  P{.\nobreak\kern 0.33333em}Patrinos, Forward-backward envelope for the sum of
  two nonconvex functions: Further properties and nonmonotone linesearch
  algorithms, \emph{SIAM Journal on Optimization} 28 (2018),  2274--2303,
  \href{https://dx.doi.org/10.1137/16M1080240}{\nolinkurl{doi:10.1137/16m1080240}}.

\bibitem{tseng2009coordinate}
P{.\nobreak\kern 0.33333em}Tseng and S{.\nobreak\kern 0.33333em}Yun, A
  coordinate gradient descent method for nonsmooth separable minimization,
  \emph{Mathematical Programming} 117 (2009),  387--423,
  \href{https://dx.doi.org/10.1007/s10107-007-0170-0}{\nolinkurl{doi:10.1007/s10107-007-0170-0}}.

\bibitem{wright2009sparse}
S.\,J{.\nobreak\kern 0.33333em}Wright, R.\,D{.\nobreak\kern 0.33333em}Nowak,
  and M.\,A.\,T{.\nobreak\kern 0.33333em}Figueiredo, Sparse reconstruction by
  separable approximation, \emph{IEEE Transactions on Signal Processing} 57
  (2009),  2479--2493,
  \href{https://dx.doi.org/10.1109/TSP.2009.2016892}{\nolinkurl{doi:10.1109/tsp.2009.2016892}}.

\bibitem{zhang2004nonmonotone}
H{.\nobreak\kern 0.33333em}Zhang and W.\,W{.\nobreak\kern 0.33333em}Hager, A
  nonmonotone line search technique and its application to unconstrained
  optimization, \emph{SIAM Journal on Optimization} 14 (2004),  1043--1056,
  \href{https://dx.doi.org/10.1137/S1052623403428208}{\nolinkurl{doi:10.1137/s1052623403428208}}.

\end{thebibliography}

\end{document}